\documentclass[11pt,a4paper]{article}
\usepackage[utf8]{inputenc}
\usepackage[T1]{fontenc}
\usepackage{lmodern}
\usepackage[english]{babel}
\usepackage{amsmath,amssymb,amsthm,mathtools}
\usepackage{enumitem}
\usepackage[a4paper,margin=2.45cm]{geometry}
\usepackage{microtype}
\usepackage[hidelinks,unicode]{hyperref}

\newtheorem{theorem}{Theorem}[section]
\newtheorem{proposition}[theorem]{Proposition}
\newtheorem{lemma}[theorem]{Lemma}
\newtheorem{corollary}[theorem]{Corollary}
\theoremstyle{remark}
\newtheorem{remark}[theorem]{Remark}

\numberwithin{equation}{section}

\newcommand{\T}{\mathbb T}
\newcommand{\Z}{\mathbb Z}

\newcommand{\Ent}{\operatorname{Ent}}
\newcommand{\supp}{\operatorname{supp}}
\newcommand{\e}{\mathrm e}
\newcommand{\1}{\mathbf 1}
\newcommand{\Ccal}{\mathcal C}

\newcommand{\dist}{\operatorname{dist}}
\newcommand{\Var}{\operatorname{Var}}

\hypersetup{
  pdftitle={Independent powers and bounded Fourier partial sums in L1},
  pdfauthor={Przemyslaw Ohrysko},
  pdfsubject={Harmonic analysis; singular measures and convolution powers},
  pdfkeywords={Fourier-Stieltjes series, Riesz products, independent powers, Rajchman measures, natural spectrum}
}

\title{Independent powers and bounded Fourier partial sums\protect\\ in $L^1$}
\author{Przemys{\l}aw Ohrysko\\[0.6em]
\small University of Warsaw\\
\small Faculty of Mathematics, Informatics and Mechanics\\
\small Institute of Mathematics\\
\small Banacha 2, 02-097 Warsaw, Poland\\[0.3em]
\small \href{mailto:P.Ohrysko@mimuw.edu.pl}{\texttt{P.Ohrysko@mimuw.edu.pl}}\\
\small ORCID: \href{https://orcid.org/0000-0002-2763-7031}{0000-0002-2763-7031}}
\date{October 25, 2026}

\begin{document}
\maketitle

\begin{abstract}
We construct a probability measure on the circle with independent convolution
powers and uniformly bounded symmetric Fourier partial sums in $L^1$.
Every positive convolution power is therefore singular with respect to Haar
measure. The measure is Rajchman. The construction
answers a question of Andersson concerning measures with bounded Fourier
partial sums outside the radical of $L^1$, and contrasts with his theorem
that the classical examples of Weiss, Katznelson, and Brown--Hewitt have
convolution squares in $L^2$. Our proof combines sparse probability blocks,
an endpoint $L\log L$ estimate, and a diagonal generalized Riesz product.
Exact factorization of local moments permits the simultaneous separation
of all convolution powers by Brown's singularity criterion. By the classical
spectral theory of independent-power measures, its spectrum is the closed
unit disc. In particular, it has non-natural spectrum and induces a
non-decomposable convolution operator on $L^1$. The image of the measure in $M(\T)/L^1(\T)$ has spectral radius one, and
the closed unital algebra it generates in $M(\T)$ is isometrically
isomorphic to $\ell^1(\mathbb N_0)$.
\end{abstract}



\section{Introduction}\label{sec:introduction}

Let $\T=\mathbb R/2\pi\mathbb Z$, equipped with normalized Haar measure $m$.
The starting point is a question of Steinhaus about a formal trigonometric
series
\[
  \sum_{n\in\Z}c_n\e^{inx}.
\]
As recorded by Helson \cite{Helson}, Steinhaus conjectured that if all
sufficiently large symmetric partial sums are nonnegative at every point,
then $c_n\to0$ as $|n|\to\infty$. Nonnegativity gives a uniform $L^1$ bound,
since the integral of every partial sum is the same constant $c_0$.
In 1954 Helson \cite{Helson} proved the stronger assertion
\begin{equation}\label{eq:Helson}
  \sup_{N\ge0}\left\|\sum_{|n|\le N}c_n\e^{inx}\right\|_1<\infty
  \quad\Longrightarrow\quad
  c_n\longrightarrow0\quad(|n|\to\infty).
\end{equation}
The boundedness assumption also ensures that the series is the
Fourier--Stieltjes series of a finite measure. Thus Helson's theorem
identifies a class of Rajchman measures by a condition on their partial
sums.

At the end of the same paper, Helson asked whether this boundedness
condition forces the representing measure to be absolutely continuous
\cite[p.~206]{Helson}. This is a stronger question than Steinhaus's original
conjecture. Weiss \cite{Weiss} answered it negatively in 1959 by constructing
a singular measure with uniformly bounded symmetric Fourier partial sums
in $L^1$. Katznelson \cite{Katznelson} subsequently obtained a singular
example for which every symmetric partial sum is positive. In particular,
even the positivity assumption in Steinhaus's question does not imply
absolute continuity. Brown and Hewitt \cite{BrownHewitt} later developed a
general construction of singular measures with positive partial sums and
prescribed admissible decay of their Fourier--Stieltjes coefficients.

For precision, let $M(\T)$ denote the convolution Banach algebra of finite
complex Borel measures on $\T$, with its total variation norm. We write
\[
  \widehat\nu(n)=\int_\T\e^{-inx}\,d\nu(x),
  \qquad
  S_N\nu(x)=\sum_{|n|\le N}\widehat\nu(n)\e^{inx},
\]
and set
\begin{equation}\label{eq:Mbdef}
  M_b(\T)=\left\{\nu\in M(\T):
  \|\nu\|_b:=\sup_{N\ge0}\|S_N\nu\|_{L^1(m)}<\infty\right\}.
\end{equation}
Here and below, by bounded partial sums we mean boundedness in the $L^1$ norm,
not uniform integrability in the measure-theoretic sense. We also use
\[
  M_0(\T)=\{\nu\in M(\T):\widehat\nu(n)\to0\text{ as }|n|\to\infty\}.
\]
Helson's theorem says that $M_b(\T)\subset M_0(\T)$, whereas the constructions
just described show that $M_b(\T)$ is not contained in $L^1(\T)$.
Quantitative refinements of Helson's theorem were obtained by Pigno and
Smith \cite{PignoSmith} and by Fournier \cite{Fournier}; the Banach-space
structure of $M_b$ and related spaces was studied by Wojtaszczyk
\cite{Wojtaszczyk} and Andersson \cite{Andersson1999}.

Andersson \cite[pp.~13--18]{Andersson} examined the constructions of Weiss,
Katznelson, and Brown--Hewitt from the viewpoint of convolution
regularization. He proved that the measures $\nu$ produced by these three
constructions satisfy
\begin{equation}\label{eq:Andersson}
  \nu*\nu\in L^2(\T).
\end{equation}
This theorem concerns those particular constructions, rather than arbitrary
members of $M_b(\T)$. To formulate his question, Andersson uses the notation
$\operatorname{Rad}L^1$. We make explicit the standard convention for the
radical of an ideal in a commutative Banach algebra; see
\cite[Symbols, Conventions, and Terminology, entry \emph{Rad I}]{GrahamMcGehee} and \cite[p.~531]{BrownMoran1983}.
Let $q:M(\T)\to M(\T)/L^1(\T)$ be the quotient map. We write
\[
  \operatorname{Rad}(L^1(\T))
  :=q^{-1}\!\left(\operatorname{rad}\bigl(M(\T)/L^1(\T)\bigr)\right),
\]
where $\operatorname{rad}$ denotes the Jacobson radical. Equivalently,
$\operatorname{Rad}(L^1(\T))$ is the intersection of all maximal ideals of
$M(\T)$ containing $L^1(\T)$, and
\[
  \nu\in\operatorname{Rad}(L^1(\T))
  \quad\Longleftrightarrow\quad
  r_{M(\T)/L^1(\T)}(q(\nu))=0.
\]
In particular, $\nu^{*k}\in L^1(\T)$ for some $k\ge1$ implies membership
in this radical, since then $q(\nu)$ is nilpotent. This sufficient
condition is distinct from the defining quasinilpotence condition.
Andersson explicitly asked whether
$M_b(\T)\setminus\operatorname{Rad}(L^1(\T))$ is nonempty
\cite[p.~18]{Andersson}. He also asked about finite regularization
thresholds beyond the convolution square.

The following result gives a measure in $M_b(\T)$ whose distinct positive
convolution powers are mutually singular. We use the standard term
\emph{independent powers} for this property; see
\cite[Chapter~6]{GrahamMcGehee} and \cite[Section~1.1]{Parreau}.
Rajchman measures with independent powers are classical: Brown's work
\cite{Brown1975} provides such Riesz products, including the example
recorded after Corollary~7.2.3 in \cite[p.~204]{GrahamMcGehee}.
The point here is to obtain independent powers together with uniformly
bounded Fourier partial sums in $L^1$.

\begin{theorem}\label{thm:main}
There exists a probability measure $\mu\in M_b(\T)\cap M_0(\T)$ with
independent convolution powers:
\[
  \mu^{*r}\perp\mu^{*s}\qquad(r,s\ge1,\ r\ne s).
\]
Consequently, every positive convolution power is singular with respect
to Haar measure:
\[
  \mu^{*r}\perp m\qquad(r\ge1).
\]
\end{theorem}

The implication in the last assertion holds for every probability measure
with independent powers on a compact group; see
Lemma~\ref{lem:independent-singular}. In particular, it does not require
the Rajchman property.

The answer to Andersson's first question holds in a stronger quotient-algebra
form. For the measure in Theorem~\ref{thm:main} and the quotient map
$q:M(\T)\to M(\T)/L^1(\T)$, we shall prove that
\begin{equation}\label{eq:introquotient}
  \|q(\mu)^r\|=1\quad(r\ge1),
  \qquad
  r_{M(\T)/L^1(\T)}(q(\mu))=1.
\end{equation}
In particular, $\mu\notin\operatorname{Rad}(L^1(\T))$, so the example
answers Andersson's nonemptiness question. The quotient computation also
quantifies the failure of regularization: the quotient powers have norm
one, rather than merely being nonzero. The present paper concerns this
nonregularizing case rather than prescribed finite regularization
thresholds.

There are also spectral and operator-theoretic consequences. A measure
$\nu\in M(\T)$ has \emph{natural spectrum} if
\[
  \sigma_{M(\T)}(\nu)=\overline{\widehat\nu(\Z)}.
\]
For $\nu\in M_0(\T)$, the set on the right is
$\widehat\nu(\Z)\cup\{0\}$. Zafran's ideal is
\begin{equation}\label{eq:zafranideal}
  \Ccal=\left\{\nu\in M_0(\T):
  \sigma_{M(\T)}(\nu)=\widehat\nu(\Z)\cup\{0\}\right\}.
\end{equation}
The fact that this set is an ideal, as well as the relevant spectral
criterion, is due to Zafran \cite{Zafran}.

\begin{corollary}\label{cor:zafran}
The measure $\mu$ in Theorem~\ref{thm:main} satisfies
\[
  \sigma_{M(\T)}(\mu)=\{z\in\mathbb C:|z|\le1\}.
\]
Consequently, it does not have natural spectrum, and
$M_b(\T)\not\subset\Ccal$.
\end{corollary}

\begin{corollary}\label{cor:decomp}
For the measure $\mu$ in Theorem~\ref{thm:main}, the convolution operator
\[
  T_\mu:L^1(\T)\longrightarrow L^1(\T),\qquad T_\mu f=f*\mu,
\]
is not decomposable in the sense of local spectral theory.
\end{corollary}

The proofs of these corollaries are given in Section~\ref{sec:spectral}.
The full-disc assertion is an application of the classical spectral theory
of independent-power measures \cite{BaileyBrownMoran}; see the precise
circle-group statement in \cite[Theorem~6.1.1(iv)]{GrahamMcGehee}.
Independent powers also give the standard isometric identification of the
closed unital subalgebra generated by $\mu$ with $\ell^1(\mathbb N_0)$
\cite[p.~138]{GrahamMcGehee}. The same is true of the operator algebra
generated by $T_\mu$ in $B(L^1(\T))$.

The proof uses several classical construction principles. The decomposition
of partial sums into complete frequency blocks and boundary pieces follows
the pattern in Brown--Hewitt \cite[Section~7.4]{BrownHewitt}; here the
boundary pieces are estimated in $L^1$ by an endpoint entropy bound rather
than by the Wiener-algebra norm. Sparse random sets with small Fourier
coefficients go back to Kaufman \cite{Kaufman}, while orthogonality-based
separation of Riesz products was developed by Peyri\`ere \cite{Peyriere}
and Brown--Moran \cite{BrownMoran1974}. The construction below combines
these principles with a diagonal choice of convolution orders and a
uniform bound for all Fourier partial sums.

The proof starts in Section~\ref{sec:analytic} with the local ingredients
for a diagonal generalized Riesz product. For each integer $r\ge1$, we
build a nonnegative trigonometric polynomial
$H$ of integral one such that $H^{*r}$ is concentrated on a set of small
Haar measure. Its Fourier coefficients at nonzero frequencies are small, with an absolute
constant independent of $r$. Diluting the block by setting
$Q=(1-\theta)+\theta H$, where $\theta$ is of order $(\log M)^{-1}$ and
$M$ bounds its degree, controls both its partial sums and its entropy.
An endpoint $L\log L$ estimate then bounds the incomplete frequency
blocks arising in a partial sum of a product.

In Section~\ref{sec:parameters}, we partition the block indices into
infinite classes $J_r$, assigning the concentration requirement of order
$r$ to the blocks in $J_r$. The scales
are chosen to give exact dissociation, including the frequencies of
polynomial tests adapted to the concentration sets. They are also chosen
to make the entropy of the partial products grow at most linearly. This
yields a weak-$*$ limit $\mu\in M_b(\T)$, and a direct coefficient estimate
shows that $\mu$ is Rajchman.

Section~\ref{sec:independent} establishes independent powers. Convolution
powers do not enlarge the Fourier support of an individual block, so the
same scales give exact factorizations of the
moments of the tests for every power. On the blocks in $J_r$, the local
perturbation in $\mu^{*r}$ has size $\theta_j^r$, whereas its size in
$\mu^{*s}$ is $\theta_j^s=o(\theta_j^r)$ when $s>r$. Brown's criterion then
separates the powers. No concentration or smoothing estimate for the
later local power $H_j^{*s}$ is required. Singularity of every positive
power with respect to Haar measure then follows from
Lemma~\ref{lem:independent-singular}.

We identify $f\in L^1(\T)$ with the measure $f\,m$, and use the same
notation for their convolution powers. In particular,
$H^{*r}(E)=\int_E H^{*r}\,dm$ for a probability density $H$.
The constant density $1$ represents Haar measure, not the identity
$\delta_0$ of $M(\T)$. We write $\|\cdot\|_1$ for the $L^1(m)$ norm,
$\mathbb N=\{1,2,\ldots\}$, and $\mathbb N_0=\{0,1,2,\ldots\}$.
All logarithms are natural, with $0\log0=0$. Constants denoted by $C$
may change from line to line; dependence on a convolution order is shown
by a subscript. The notation $u_j\asymp_r v_j$ means two-sided estimates
with positive constants depending only on the fixed order $r$.

\section{Analytic estimates and sparse blocks}\label{sec:analytic}
\label{sec:sparse}
We begin with two estimates for the partial sums of a product and then
construct the local probability blocks. For a trigonometric polynomial
$f$ and an integer interval $I\subset\Z$, put
\[
  P_I f=\sum_{n\in I}\widehat f(n)\e^{inx}.
\]
For a nonnegative polynomial $R$ with $\int R\,dm=1$, define
\[
  \Ent(R)=\int_\T R\log R\,dm.
\]
By convexity, $\Ent(R)\ge0$.

\begin{lemma}[Endpoint estimate]\label{lem:LlogL}
There is an absolute constant $C>0$ such that
\[
  \|P_I R\|_1\le C\bigl(1+\Ent(R)\bigr)
\]
for every nonnegative trigonometric polynomial $R$ with integral one and
every integer interval $I\subset\Z$.
\end{lemma}

\begin{proof}
We use the following quantitative consequence of Zygmund's
conjugate-function estimate \cite[pp.~285--286, 291--294]{Zygmund1929}:
\begin{equation}\label{eq:Zygmundnormalized}
  \|\widetilde f\|_1\le
  C_0\left(1+\int_\T |f|\log^+|f|\,dm\right),
  \qquad \|f\|_1=1.
\end{equation}
The estimate holds for complex trigonometric polynomials as well: apply
its real-variable form after normalizing the real and imaginary parts
separately (the zero-norm case being trivial). This changes only the
absolute constant. The identity expressing the Riesz projection $P_+$ in terms of the conjugate
function yields the same bound for $P_+f$, where $P_+$ retains the
nonnegative frequencies.
For a finite interval $I=[a,b]\cap\Z$, we have
\[
  P_I f=\e^{iax}P_+(\e^{-iax}f)
       -\e^{i(b+1)x}P_+(\e^{-i(b+1)x}f).
\]
Modulation preserves both $|f|$ and $\|f\|_1$. Applying this identity to
$R$ and observing that
\[
  \int R\log^+R\,dm
  =\Ent(R)+\int_{\{R<1\}}(-R\log R)\,dm
  \le\Ent(R)+\e^{-1}
\]
proves the assertion. For an unbounded interval, only the finite
intersection with the Fourier support of $R$ matters.
\end{proof}

The following estimate uses the same complete-block and boundary-block
decomposition as \cite[Section~7.4, especially (7.4.3)]{BrownHewitt}.
The endpoint estimate above replaces the pointwise control of the
boundary term used there by an $L^1$ bound in terms of entropy.

\begin{lemma}[Separated blocks]\label{lem:block}
Let
\[
  R(x)=\sum_{|n|\le D}r_n\e^{inx}\ge0,
  \qquad \int R\,dm=1,
\]
and let
\[
  Q(x)=\sum_{|j|\le M}q_j\e^{ijx},\qquad q_0=1.
\]
If $\lambda\in\mathbb N$ and $\lambda>2D$, then, with
$a(Q)=\max_{j\ne0}|q_j|$,
\begin{equation}\label{eq:blockbound}
  \|R\,Q(\lambda\cdot)\|_b
  \le\max\left\{\|R\|_b,
    \|Q\|_b+2C\,a(Q)(1+\Ent R)\right\},
\end{equation}
where $C$ is the constant in Lemma~\ref{lem:LlogL}.
\end{lemma}

\begin{proof}
The Fourier support of the product is contained in the disjoint blocks
\[
  j\lambda+[-D,D],\qquad |j|\le M.
\]
If $N<D$, then
\[
  [-N,N]\cap\bigl(j\lambda+[-D,D]\bigr)=\varnothing
  \qquad(j\ne0),
\]
so that $S_N(RQ(\lambda\cdot))=S_NR$. Its $L^1$ norm is therefore at most
$\|R\|_b$.

Suppose now that $N\ge D$. The central block is then complete, and a
symmetric truncation at $[-N,N]$ cuts at most two boundary blocks, each
with a nonzero index $j$. The complete blocks contribute
$R(x)T(\lambda x)$, where $T$ is a symmetric partial sum of $Q$.
The function $|T(\lambda\cdot)|$ has Fourier coefficients only at
frequencies in $\lambda\Z$, whereas $R$ has no nonzero frequency in
$\lambda\Z$. Hence
\[
  \int_\T R(x)|T(\lambda x)|\,dm(x)=\int_\T|T|\,dm\le\|Q\|_b.
\]
This identity follows directly by integrating against the finite Fourier
expansion of $R$; no pointwise convergence of a Fourier series is needed.
Each incomplete boundary block is of the form
\[
  q_j\e^{ij\lambda x}P_I R(x)
\]
for some interval $I\subset[-D,D]$ and some $j\ne0$.
Lemma~\ref{lem:LlogL} bounds its norm by
$C\,a(Q)(1+\Ent R)$. Combining the two cases proves
\eqref{eq:blockbound}.
\end{proof}

We now construct the local blocks to which these estimates will be
applied. Each block must have small Fourier coefficients at nonzero frequencies while
its prescribed convolution power is concentrated on a set of small Haar
measure. No regularity of the next power is required.
The probabilistic small-support/small-Fourier-coefficient principle has
an antecedent in Kaufman's finite-group construction \cite{Kaufman}.
For a related use of random discrete measures followed by smoothing,
with the different aim of making convolution squares regular, see
K\"orner \cite[Sections~2--3]{Korner2008}. Here Hoeffding's inequality
and Fej\'er smoothing give the precise bounds needed at an arbitrary
prescribed convolution order.

\begin{lemma}\label{lem:sparse}
There is an absolute constant $A>0$ such that, for every integer $r\ge1$,
there are constants $C_r>0$ and $s_0(r)\ge2$ with the following property.
For every integer $s\ge s_0(r)$, set
\[
  M=\left\lceil s^{r+1/2}\right\rceil.
\]
There exist a nonnegative trigonometric polynomial $H$ of degree less than
$M$, with integral one, and a Borel set $E\subset\T$ such that
\begin{align}
  0\le H&\le M,\label{eq:Hpoint}\\
  a(H):=\max_{n\ne0}|\widehat H(n)|
     &\le A\sqrt{\frac{\log M}{s}},\label{eq:Hcoeff}\\
  H^{*r}(E)&\ge\frac78,\label{eq:Hmass}\\
  m(E)&\le C_rs^{-1/2}.\label{eq:Esize}
\end{align}
\end{lemma}

\begin{proof}
Let $X_1,\ldots,X_s$ be independent Haar-distributed random variables on
$\T$, and put
\[
  \nu_\omega=\frac1s\sum_{\ell=1}^s\delta_{X_\ell(\omega)}.
\]
For fixed $n\ne0$, the random variables $\e^{-inX_\ell}$ are independent,
have modulus one, and have mean zero. Hoeffding's inequality
\cite{Hoeffding}, applied to the real and imaginary parts, gives
\[
  \mathbb P\bigl(|\widehat\nu_\omega(n)|>t\bigr)
  \le4\exp(-c_0st^2)
\]
with an absolute constant $c_0>0$. Taking a union bound over
$1\le|n|<M$ and setting $t=A\sqrt{(\log M)/s}$ yields
\[
  \mathbb P\left(\max_{1\le|n|<M}|\widehat\nu_\omega(n)|>t\right)
  \le8M^{1-c_0A^2}.
\]
We choose the absolute constant $A$ so that $c_0A^2>2$ and then enlarge
$s_0(r)$ if necessary. The right-hand side is less than one. Thus there
is a realization
\[
  \nu=\frac1s\sum_{\ell=1}^s\delta_{x_\ell}
\]
for which all these coefficients have modulus at most $t$.

Let $F_{M-1}$ be the normalized Fej\'er kernel of degree $M-1$,
\[
  F_{M-1}(x)=\frac1M\left|\sum_{k=0}^{M-1}\e^{ikx}\right|^2,
\]
and define
\[
  H=\nu*F_{M-1}.
\]
Then $H\ge0$, $\int H\,dm=1$, $\deg H<M$, and
\[
  H(x)=\frac1s\sum_{\ell=1}^sF_{M-1}(x-x_\ell),
  \qquad 0\le H(x)\le M.
\]
Moreover,
\[
  \widehat H(n)=
  \begin{cases}
    (1-|n|/M)\widehat\nu(n),& |n|<M,\\
    0,& |n|\ge M,
  \end{cases}
\]
which proves \eqref{eq:Hcoeff}.

It remains to establish concentration. We have
\[
  H^{*r}=\nu^{*r}*F_{M-1}^{*r},
\]
and $\nu^{*r}$ is supported on
\[
  Y=\{x_{i_1}+\cdots+x_{i_r}:1\le i_1,\ldots,i_r\le s\},
  \qquad |Y|\le s^r.
\]
For representatives $x\in[-\pi,\pi]$, the Fej\'er kernel satisfies
\[
  F_{M-1}(x)\le C\min\left\{M,\frac1{Mx^2}\right\}.
\]
Consequently, for each $\eta>0$ there is $L=L(\eta)$, independent of $M$,
such that
\[
  \int_{|x|\le L/M}F_{M-1}(x)\,dm(x)\ge1-\eta.
\]
Choose $\eta=(8r)^{-1}$ and enlarge $s_0(r)$ so that $rL/M<\pi$.
If $Y_1,\ldots,Y_r$ are independent with density $F_{M-1}$, a union
bound shows that, with probability at least $7/8$, all their
representatives satisfy $|Y_q|\le L/M$. On this event the sum belongs to
the arc $I_r=[-rL/M,rL/M]$. Therefore
\[
  F_{M-1}^{*r}(I_r)\ge\frac78,
  \qquad m(I_r)\le\frac{C_r}{M}.
\]
Define $E=\bigcup_{y\in Y}(y+I_r)$. Then
\[
  m(E)\le C_r\frac{s^r}{M}\le C_rs^{-1/2},
\]
while $E-y$ contains $I_r$ for each $y\in Y$. It follows that
\[
  H^{*r}(E)=\int F_{M-1}^{*r}(E-y)\,d\nu^{*r}(y)\ge\frac78.
\]
The constant in the coefficient estimate is independent of $r$;
only $C_r$ and the threshold $s_0(r)$ depend on the order.
\end{proof}

\begin{remark}[A symmetric version]\label{rem:symmetric}
The construction may be made symmetric without changing any of the later
arguments. In the preceding proof, replace the empirical measure by
\[
  \nu_\omega^{\rm sym}
  =\frac1{2s}\sum_{\ell=1}^s
    (\delta_{X_\ell(\omega)}+\delta_{-X_\ell(\omega)}).
\]
Then
\[
  \widehat{\nu_\omega^{\rm sym}}(n)
  =\frac1s\sum_{\ell=1}^s\cos(nX_\ell(\omega))
\]
is real, and Hoeffding's inequality gives the same coefficient bound.
The support of $(\nu_\omega^{\rm sym})^{*r}$ has cardinality at most
$(2s)^r$, so the factor $2^r$ is absorbed into $C_r$ in
\eqref{eq:Esize}. The Fej\'er smoothing is symmetric as well. Thus all
local blocks $H_j$, hence all $Q_j$ and $R_N$, may be chosen symmetric,
and the weak-$*$ limit satisfies $\widetilde\mu=\mu$. Consequently,
Theorem~\ref{thm:main} remains valid with the additional requirement that
$\mu$ be symmetric (equivalently, $\widehat\mu(n)\in\mathbb R$ for all
$n\in\mathbb Z$).
\end{remark}

\section{Construction of the measure}\label{sec:parameters}
\label{sec:scales}\label{sec:entropy}\label{sec:limit}
We combine the blocks from Lemma~\ref{lem:sparse} into a single generalized
Riesz product. The block indices are assigned to convolution orders, and
the dilations are chosen to ensure both exact dissociation and uniform
control of the partial sums. The dissociation and coefficient-stabilization
arguments are standard for generalized Riesz products; compare
\cite[Sections~7.3 and~8.1--8.2]{BrownHewitt} and
\cite[Section~2.3]{Parreau}. We retain the details because the frequency
bounds here must also accommodate all the polynomial tests.

Partition the positive integers into the disjoint infinite classes
\[
  J_r=\{j\ge1:v_2(j)=r-1\}
      =\{2^{r-1}(2m-1):m\ge1\},\qquad r\ge1,
\]
where $v_2(j)$ is the exponent of $2$ in the factorization of $j$. Put
$r_j=1+v_2(j)$, so that $j\in J_{r_j}$.

Fix $c\in(0,1/10)$. For each $r\ge1$, choose an integer $S_r\ge s_0(r)$
so large that, for every integer $s\ge S_r$, the quantities
\[
  M_r(s)=\left\lceil s^{r+1/2}\right\rceil,
  \qquad
  \vartheta_r(s)=\frac{c}{\log M_r(s)}
\]
satisfy
\begin{align}
  0<\vartheta_r(s)&\le\frac12,\label{eq:S1}\\
  C_rs^{-1/2}&\le\vartheta_r(s)^{2r}.\label{eq:S2}
\end{align}
This is possible because, for each fixed $r$,
$s^{-1/2}=o((\log s)^{-2r})$ as $s\to\infty$.
For $j\ge1$, set
\[
  s_j=\max\{j^4,S_{r_j}\},
  \qquad M_j=\left\lceil s_j^{r_j+1/2}\right\rceil.
\]
Choose a block $H_j$ and a set $E_j$ from Lemma~\ref{lem:sparse}, with
order $r_j$ and parameter $s_j$, and define
\begin{equation}\label{eq:Qdef}
  \theta_j=\frac{c}{\log M_j},
  \qquad Q_j=(1-\theta_j)+\theta_jH_j.
\end{equation}

\begin{lemma}[Diluted blocks]\label{lem:Q}
There is a constant $C>0$, independent of $j$ and $r_j$, such that
\begin{align}
  Q_j&>0,\qquad \int Q_j\,dm=1,\label{eq:Qprob}\\
  \|Q_j\|_b&\le C,\label{eq:Qb}\\
  a(Q_j)&\le\frac{C}{\sqrt{s_j\log M_j}}
          \le\frac{C}{j^2},\label{eq:aQ}\\
  \Ent(Q_j)&\le c.\label{eq:EntQ}
\end{align}
Moreover, for every integer $\ell\ge1$,
\begin{equation}\label{eq:Qpower}
  Q_j^{*\ell}=(1-\theta_j^\ell)+\theta_j^\ell H_j^{*\ell}.
\end{equation}
\end{lemma}

\begin{proof}
The first assertion follows from \eqref{eq:S1} and the normalization of
$H_j$. Let $D_N=\sum_{|n|\le N}\e^{inx}$ be the Dirichlet kernel.
For $N<M_j$, Young's inequality gives
\[
  \|S_NH_j\|_1=\|H_j*D_N\|_1
  \le\|D_N\|_1\le C\log(M_j+1).
\]
For $N\ge M_j-1$, the partial sum is $H_j$ itself and has norm one.
Consequently
\[
  \|Q_j\|_b\le 1+C\theta_j\log M_j\le1+Cc,
\]
which proves \eqref{eq:Qb}. The absolute constant in
\eqref{eq:Hcoeff} yields
\[
  a(Q_j)=\theta_j a(H_j)
  \le\frac{Ac}{\sqrt{s_j\log M_j}}.
\]
Since $s_j\ge j^4$, this proves \eqref{eq:aQ}, after a change of constant.

By convexity of $\Phi(t)=t\log t$ and the identity $\Phi(1)=0$,
\[
  \Phi\bigl((1-\theta_j)\cdot1+\theta_jH_j\bigr)
  \le\theta_j\Phi(H_j).
\]
The bounds $0\le H_j\le M_j$ and $\int H_j\,dm=1$ therefore imply
\[
  \Ent(Q_j)\le\theta_j\Ent(H_j)
  \le\theta_j\log M_j=c.
\]
Finally, for $n\ne0$ we have
$\widehat Q_j(n)=\theta_j\widehat H_j(n)$. Raising this identity to the
power $\ell$ proves equality of the Fourier coefficients at nonzero frequencies in
\eqref{eq:Qpower}; both Fourier coefficients at frequency zero equal one.
\end{proof}

Set $\delta_j=\theta_j^{r_j}$. By \eqref{eq:S2} and
Lemma~\ref{lem:sparse},
\begin{equation}\label{eq:concentrationsets}
  H_j^{*r_j}(E_j)\ge\frac78,
  \qquad m(E_j)\le\delta_j^2.
\end{equation}
We shall use polynomial tests rather than the indicators of these sets.
Choose $T_j\in\mathbb N$ sufficiently large and put
\begin{equation}\label{eq:pj}
  p_j=\1_{E_j}*F_{T_j}.
\end{equation}
Then $p_j$ is a real nonnegative trigonometric polynomial of degree at
most $T_j$, and $T_j$ can be chosen so that
\begin{align}
  0\le p_j\le1,\qquad
  a_j:=\int p_j\,dm=m(E_j)&\le\delta_j^2,\label{eq:aj}\\
  b_j:=\int p_jH_j^{*r_j}\,dm&\ge\frac34.\label{eq:bj}
\end{align}
Indeed, the Fej\'er kernels form an approximate identity and
$H_j^{*r_j}$ is continuous. Hence
\[
  \int (\1_{E_j}*F_T)H_j^{*r_j}\,dm
  \longrightarrow\int\1_{E_j}H_j^{*r_j}\,dm\ge\frac78
  \qquad(T\to\infty).
\]
This proves that the required choice of $T_j$ is possible.

\begin{remark}\label{rem:thresholds}
The separate thresholds $S_r$ ensure
$m(E_j)\le\theta_j^{2r_j}$ for every block, even though the orders $r_j$
are unbounded. A single assertion valid for sufficiently large $j$ at a
fixed order would not suffice. On the other hand, on any fixed class
$J_r$, we have $s_j=j^4$ for all sufficiently large $j$. This will preserve
the divergence needed to separate distinct convolution powers.
\end{remark}

Having fixed the blocks and their polynomial tests, we choose their
frequency scales. The tests must be included in the dissociation
condition, so we define the frequency bounds
\begin{equation}\label{eq:budget}
  B_j=2(M_j+T_j+1).
\end{equation}
We choose integers $\lambda_j\ge1$ inductively, requiring
\begin{equation}\label{eq:separation}
  \lambda_j>4\sum_{i<j}B_i\lambda_i.
\end{equation}
At each step, the scale will also be chosen sufficiently large to satisfy
the entropy estimate in Lemma~\ref{lem:entropy} below. Both requirements
are lower bounds on the next scale and are therefore compatible.

\begin{lemma}[Exact dissociation]\label{lem:diss}
If
\[
  \sum_{i=1}^N c_i\lambda_i=0,
  \qquad c_i\in\Z,\qquad |c_i|\le B_i,
\]
then $c_i=0$ for every $i$.
\end{lemma}

\begin{proof}
If some coefficient is nonzero, let $q$ be its largest index. Since
$c_q$ is a nonzero integer, \eqref{eq:separation} gives
\[
  |c_q|\lambda_q\ge\lambda_q
  >4\sum_{i<q}B_i\lambda_i
  \ge\left|\sum_{i<q}c_i\lambda_i\right|,
\]
contrary to the vanishing of the sum.
\end{proof}

In particular, any finite representation
\begin{equation}\label{eq:uniquerepresentation}
  n=\sum_i c_i\lambda_i,
  \qquad c_i\in\Z,\qquad |c_i|\le M_i+T_i+1,
\end{equation}
is unique, if it exists. To see this, subtract two such representations
and apply Lemma~\ref{lem:diss} with the bounds \eqref{eq:budget}.

To control the partial sums of the full product, we combine this
dissociation with an entropy estimate. Put $R_0=1$ and
\begin{equation}\label{eq:products}
  R_N(x)=\prod_{j=1}^N Q_j(\lambda_jx),\qquad N\ge1.
\end{equation}
Each $R_N$ is strictly positive. Exact dissociation shows that its zero
Fourier coefficient is the product of the Fourier coefficients at frequency zero of the
$Q_j$, so that $\int R_N\,dm=1$.

\begin{lemma}[Entropy growth]\label{lem:entropy}
The scales can be chosen to satisfy \eqref{eq:separation} and
\begin{equation}\label{eq:entropy}
  \Ent(R_N)\le cN+1\qquad(N\ge0).
\end{equation}
\end{lemma}

\begin{proof}
Suppose that $\lambda_1,\ldots,\lambda_{j-1}$ have been chosen.
The function $R_{j-1}\log R_{j-1}$ is continuous. Since $Q_j$ is a
fixed trigonometric polynomial of integral one, the Riemann--Lebesgue
lemma implies
\[
  \int R_{j-1}\log R_{j-1}\,Q_j(\lambda\cdot)\,dm
  \longrightarrow\Ent(R_{j-1})
\]
as the positive integer $\lambda$ tends to infinity. We may therefore
choose $\lambda_j$ satisfying \eqref{eq:separation} so that the absolute
value of the difference is less than $2^{-j}$.

Also, $\lambda_j>\deg R_{j-1}$. The continuous function
$Q_j(\lambda_j\cdot)\log Q_j(\lambda_j\cdot)$ has Fourier coefficients
only at multiples of $\lambda_j$. Pairing it with the polynomial
$R_{j-1}$ therefore gives the exact identity
\[
  \int R_{j-1}Q_j(\lambda_j\cdot)
           \log Q_j(\lambda_j\cdot)\,dm=\Ent(Q_j).
\]
Expanding $R_j\log R_j$ and using \eqref{eq:EntQ}, we obtain
\[
  \Ent(R_j)
  \le\Ent(R_{j-1})+\Ent(Q_j)+2^{-j}
  \le\Ent(R_{j-1})+c+2^{-j}.
\]
Since $\Ent(R_0)=0$, summation proves \eqref{eq:entropy}.
\end{proof}

The linear entropy bound, together with the decay of $a(Q_j)$, now gives
a uniform estimate for all partial products.

\begin{proposition}[Uniformly bounded partial sums]\label{prop:Mb}
There is a constant $C>0$ such that
\begin{equation}\label{eq:uniform}
  \sup_{N\ge0}\|R_N\|_b\le C.
\end{equation}
\end{proposition}

\begin{proof}
Apply Lemma~\ref{lem:block} to $R_{j-1}$ and $Q_j$. Condition
\eqref{eq:separation} implies $\lambda_j>2\deg R_{j-1}$. By
Lemmas~\ref{lem:Q} and \ref{lem:entropy},
\[
  a(Q_j)(1+\Ent R_{j-1})
  \le\frac{C(1+j)}{j^2}\le C,
\]
while $\|Q_j\|_b\le C$. Thus
\[
  \|R_j\|_b\le\max\{\|R_{j-1}\|_b,C\},
\]
with the same constant for every $j$. Starting with $R_0=1$ proves the
assertion. In particular, the errors from the boundary blocks do not
accumulate over the scales.
\end{proof}

We pass from the partial products to the desired probability measure.
For each fixed $n\in\Z$, the coefficient $\widehat R_N(n)$ is eventually
constant. Indeed, a frequency occurring in $R_N$ has a representation
$n=\sum_{i\le N}c_i\lambda_i$ with $c_i\in\Z$ and $|c_i|\le M_i$,
which is unique by Lemma~\ref{lem:diss}. If $n\ne0$ and $q$ is the
largest index for which $c_q\ne0$, then
\[
  |n|\ge\lambda_q-\sum_{i<q}M_i\lambda_i
  >\frac34\lambda_q.
\]
Since $\lambda_q\to\infty$, no sufficiently late scale can occur in a
representation of the fixed frequency $n$. The zeroth Fourier coefficient is one
for every $N$.

The probability measures $R_Nm$ lie in a weak-$*$ compact subset of
$M(\T)$. Their cluster points have the same Fourier--Stieltjes
coefficients and hence coincide. Consequently the entire sequence
converges:
\begin{equation}\label{eq:weaklimit}
  R_Nm\stackrel{w^*}{\longrightarrow}\mu
\end{equation}
for a probability measure $\mu$. For each fixed $L$, all coefficients
with $|n|\le L$ have stabilized once $N$ is sufficiently large. Hence
\[
  S_L\mu=S_LR_N
\]
for all such $N$. Proposition~\ref{prop:Mb} gives
\begin{equation}\label{eq:muMb}
  \mu\in M_b(\T).
\end{equation}

\begin{lemma}\label{lem:Rajchman}
The measure $\mu$ belongs to $M_0(\T)$.
\end{lemma}

\begin{proof}
Suppose that $n\ne0$ and $\widehat\mu(n)\ne0$. Then $n$ has a unique
finite representation
\[
  n=\sum_{i\le j}c_i\lambda_i,
  \qquad c_j\ne0,\qquad |c_i|\le M_i,
\]
and
\[
  \widehat\mu(n)=\prod_{i\le j}\widehat Q_i(c_i).
\]
Each $Q_i$ is a probability density, so
$|\widehat Q_i(c_i)|\le1$. At the largest active index,
$|\widehat Q_j(c_j)|\le a(Q_j)$. Therefore
\begin{equation}\label{eq:Rajchmanbound}
  |\widehat\mu(n)|\le a(Q_j)\le Cj^{-2}.
\end{equation}
As $|n|\to\infty$ through frequencies with nonzero coefficient, the
largest active index must tend to infinity: bounded indices allow only
finitely many representations with the stated coefficient bounds.
This proves $\widehat\mu(n)\to0$.
\end{proof}

Although membership in $M_0(\T)$ also follows from Helson's theorem,
\eqref{eq:Rajchmanbound} proves it directly from the construction.

\section{Independent convolution powers}\label{sec:independent}
\label{sec:factorization}\label{sec:singularity}
It remains to prove that the measure constructed in
Section~\ref{sec:parameters} has independent powers. We first factor the
convolution powers and the moments of the polynomial tests. These
identities will allow us to apply Brown's criterion to any pair of
distinct powers.

For $r\ge1$ and $i\ge1$, write $G_i^{(r)}=Q_i^{*r}$. Formula
\eqref{eq:Qpower} gives
\begin{equation}\label{eq:Gdef}
  G_i^{(r)}=(1-\theta_i^r)+\theta_i^rH_i^{*r},
  \qquad \supp\widehat{G_i^{(r)}}\subset[-M_i,M_i].
\end{equation}
In particular, taking a convolution power does not enlarge the Fourier
support of a local factor. The next identity is the finite-product form
of the standard convolution rule for generalized Riesz products; see
\cite[Lemma~2.2(a)]{Parreau}.

\begin{lemma}[Factorization of powers]\label{lem:powerfactor}
For every $N,r\ge1$,
\begin{equation}\label{eq:powerfactor}
  R_N^{*r}(x)=\prod_{i=1}^N G_i^{(r)}(\lambda_i x).
\end{equation}
Moreover,
\begin{equation}\label{eq:powerlimit}
  R_N^{*r}m\stackrel{w^*}{\longrightarrow}\mu^{*r}.
\end{equation}
\end{lemma}

\begin{proof}
Every frequency occurring in $R_N$ has at most one representation
$n=\sum_{i=1}^N c_i\lambda_i$ with $|c_i|\le M_i$. For a frequency
with such a representation,
\[
  \widehat{R_N^{*r}}(n)
  =\widehat R_N(n)^r
  =\prod_{i=1}^N\widehat Q_i(c_i)^r
  =\prod_{i=1}^N\widehat{G_i^{(r)}}(c_i).
\]
By the same dissociation, the last product is the Fourier coefficient of
the right-hand side of \eqref{eq:powerfactor}. If no representation
exists, both coefficients vanish. Thus \eqref{eq:powerfactor} holds.
For every $n\in\Z$,
\[
  \widehat{R_N^{*r}}(n)=\widehat R_N(n)^r
  \longrightarrow\widehat\mu(n)^r=\widehat{\mu^{*r}}(n).
\]
Since the measures $R_N^{*r}m$ are probabilities, weak-$*$ compactness and
uniqueness of the Fourier--Stieltjes transform yield
\eqref{eq:powerlimit}.
\end{proof}

We next record exact moment identities for the polynomial tests chosen
in \eqref{eq:pj}. The use of dissociation to factor these moments and
annihilate off-diagonal covariances is the same general mechanism as in
\cite[proof of Theorem~8.3]{BrownHewitt}; our frequency bounds include
the additional tests and apply at every convolution order.

\begin{lemma}[Factorization of moments]\label{lem:moments}
Fix an integer $r\ge1$. For every $j$,
\begin{equation}\label{eq:firstmoment}
  \int p_j(\lambda_jx)\,d\mu^{*r}(x)
  =\int p_jG_j^{(r)}\,dm.
\end{equation}
For $i\ne j$,
\begin{equation}\label{eq:crossmoment}
  \int p_i(\lambda_ix)p_j(\lambda_jx)\,d\mu^{*r}(x)
  =\left(\int p_iG_i^{(r)}\,dm\right)
   \left(\int p_jG_j^{(r)}\,dm\right).
\end{equation}
The corresponding identities hold for $m$, with every $G_i^{(r)}$
replaced by $1$.
\end{lemma}

\begin{proof}
First replace $\mu^{*r}$ by $R_N^{*r}m$, where $N$ is at least as large
as all test indices. Use \eqref{eq:powerfactor}. In the first moment,
the local factor at scale $j$ is $p_jG_j^{(r)}$, whose Fourier support
is contained in
\[
  [-(M_j+T_j),M_j+T_j].
\]
All other local factors have Fourier support in $[-M_i,M_i]$.
These bounds fall within \eqref{eq:budget}. Lemma~\ref{lem:diss}
therefore implies that the zeroth Fourier coefficient of the full product can
only arise by choosing the zeroth Fourier coefficient of each local factor.
Since $\int G_i^{(r)}\,dm=1$, this proves \eqref{eq:firstmoment} for
$R_N^{*r}m$. The tests are continuous, so \eqref{eq:powerlimit} permits
passage to $\mu^{*r}$.

For \eqref{eq:crossmoment}, use the two local factors
$p_iG_i^{(r)}$ and $p_jG_j^{(r)}$ and the same argument. Setting all the
$G_i^{(r)}$ equal to $1$ proves the statements for $m$.
\end{proof}

\begin{remark}
The same scales work for every convolution order. This is because
$\widehat{Q_i^{*r}}(n)=\widehat Q_i(n)^r$, so no new local frequencies
appear. Also, when a test from $J_r$ is integrated against $\mu^{*r}$,
all untested factors contribute only their zeroth Fourier coefficient, regardless
of their assigned orders. Thus the identities are exact, not asymptotic.
\end{remark}

The moment identities reduce the separation of powers to a comparison
of their local means. The following divergence estimate will ensure that
the squared differences of the normalized means have an infinite sum.

\begin{lemma}\label{lem:divergence}
For every fixed integer $r\ge1$,
\begin{equation}\label{eq:divergence}
  \sum_{j\in J_r}\theta_j^r=\infty.
\end{equation}
\end{lemma}

\begin{proof}
Since $S_r$ is fixed, for all sufficiently large $j\in J_r$ we have
$s_j=j^4$ and
\[
  M_j=\left\lceil j^{4r+2}\right\rceil.
\]
Consequently $\theta_j^r\asymp_r(\log j)^{-r}$. The arithmetic
progression $J_r=\{2^{r-1}(2m-1):m\ge1\}$ gives the same divergence
as the series $\sum_{m\ge2}(\log m)^{-r}$, whose terms eventually
exceed $1/m$.
\end{proof}

We use Brown's singularity criterion
\cite[Section~3, Theorem~I]{Brown1977}, also recorded by Brown and
Hewitt \cite[Theorem~4.2]{BrownHewitt}. Related orthogonality arguments
for ordinary Riesz products appear in
\cite[Theorem~1.2]{Peyriere} and \cite{BrownMoran1974}.
We state the complex-valued version, with conjugates in the covariance,
and include its short proof for completeness.

\begin{theorem}[Brown's criterion]\label{thm:Brown}
Let $\nu_1,\nu_2$ be probability measures on the same measurable space
$(X,\Sigma)$, and let $\phi_n\in L^2(\nu_1)\cap L^2(\nu_2)$ for $n\ge1$. Suppose
that
\begin{equation}\label{eq:Browncov}
  \sum_{q=0}^\infty
  \sup_{\substack{n\ge1\\ k=1,2}}
  \left|
    \int\phi_n\overline{\phi_{n+q}}\,d\nu_k
    -\left(\int\phi_n\,d\nu_k\right)
     \overline{\left(\int\phi_{n+q}\,d\nu_k\right)}
  \right|<\infty
\end{equation}
and
\begin{equation}\label{eq:Brownmeans}
  \sum_{n=1}^\infty
  \left|\int\phi_n\,d\nu_1-\int\phi_n\,d\nu_2\right|^2=\infty.
\end{equation}
Then $\nu_1\perp\nu_2$.
\end{theorem}

\begin{proof}
Let $c_q$ denote the supremum in \eqref{eq:Browncov}, and put
$K=c_0+2\sum_{q\ge1}c_q<\infty$. For a finite complex scalar sequence
$(t_n)$ and $F=\sum_n t_n\phi_n$, expansion of the variance gives
\[
  \Var_{\nu_k}(F)
  \le c_0\sum_n|t_n|^2
      +2\sum_{q\ge1}c_q\sum_n|t_n||t_{n+q}|
  \le K\sum_n|t_n|^2,
\]
where the coefficients outside the finite sequence are taken to be zero.
Here $\Var_\nu(F)=\int|F-\int F\,d\nu|^2\,d\nu$.

If the measures were not mutually singular, there would exist a nonzero
positive measure $\rho$ with $\rho\le\nu_1$ and $\rho\le\nu_2$.
Write $\alpha=\rho(X)>0$. By Cauchy--Schwarz,
\[
  \left|\int F\,d\nu_k-\frac1\alpha\int F\,d\rho\right|
  \le\alpha^{-1/2}\bigl(\Var_{\nu_k}(F)\bigr)^{1/2}.
\]
Subtracting the two means, we obtain
\[
  \left|\sum_n t_n
    \left(\int\phi_n\,d\nu_1-\int\phi_n\,d\nu_2\right)\right|
  \le2\sqrt{\frac K\alpha}\left(\sum_n|t_n|^2\right)^{1/2}.
\]
By $\ell^2$ duality, the sequence of differences of the means belongs
to $\ell^2$, contradicting \eqref{eq:Brownmeans}.
\end{proof}

We apply Brown's criterion to the constructed measure. To compare two
powers $r<s$, it suffices to use the blocks in $J_r$.
Only positivity and normalization will be required of $H_j^{*s}$;
no concentration estimate at the later order is needed.

\begin{proposition}[Independent powers]\label{prop:independent}
For all integers $1\le r<s$,
\[
  \mu^{*r}\perp\mu^{*s}.
\]
\end{proposition}

\begin{proof}
Fix $r<s$. For $j\in J_r$, put
\[
  \delta_j=\theta_j^r,\qquad \varepsilon_j=\theta_j^s.
\]
Then
\begin{equation}\label{eq:epsdelta}
  \frac{\varepsilon_j}{\delta_j}=\theta_j^{s-r}\longrightarrow0
  \qquad(j\to\infty,\ j\in J_r),
\end{equation}
whereas Lemma~\ref{lem:divergence} gives
\begin{equation}\label{eq:deltadiv-ind}
  \sum_{j\in J_r}\delta_j=\infty.
\end{equation}
Retain the test $p_j$ from \eqref{eq:pj} and recall that
\begin{equation}\label{eq:ab-ind}
  a_j=\int p_j\,dm\le\delta_j^2,
  \qquad b_j=\int p_jH_j^{*r}\,dm\ge\frac34.
\end{equation}
For the later power, define only
\begin{equation}\label{eq:c-ind}
  c_j=\int p_jH_j^{*s}\,dm\in[0,1].
\end{equation}
The indicated bounds follow from $0\le p_j\le1$ and the fact that
$H_j^{*s}$ is a probability density.

By Lemma~\ref{lem:moments} and \eqref{eq:Gdef},
\begin{align}
  \int p_j(\lambda_jx)\,d\mu^{*r}(x)
  &=(1-\delta_j)a_j+\delta_jb_j,\label{eq:mean-r-ind}\\
  \int p_j(\lambda_jx)\,d\mu^{*s}(x)
  &=(1-\varepsilon_j)a_j+\varepsilon_jc_j.\label{eq:mean-s-ind}
\end{align}
Their difference is
$\delta_j(b_j-a_j)-\varepsilon_j(c_j-a_j)$. Since
$|c_j-a_j|\le1$, it follows that
\begin{align}
  &\left|\int p_j(\lambda_jx)\,d\mu^{*r}
        -\int p_j(\lambda_jx)\,d\mu^{*s}\right|\notag\\
  &\hspace{2em}\ge
    \delta_j(b_j-a_j)-\varepsilon_j|c_j-a_j|
    \ge\delta_j(b_j-a_j)-\varepsilon_j.\label{eq:mean-diff-ind}
\end{align}
By \eqref{eq:epsdelta} and \eqref{eq:ab-ind}, after discarding finitely
many indices in $J_r$ we may assume that
\[
  b_j-a_j\ge\frac23,\qquad
  \varepsilon_j\le\frac16\delta_j.
\]
Consequently
\begin{equation}\label{eq:mean-lower-ind}
  \left|\int p_j(\lambda_jx)\,d\mu^{*r}
        -\int p_j(\lambda_jx)\,d\mu^{*s}\right|
  \ge\frac12\delta_j.
\end{equation}
Enumerate the remaining tail as $J_r'=\{j_1<j_2<\cdots\}$, and set
\begin{equation}\label{eq:phi-ind}
  \phi_n(x)=\delta_{j_n}^{-1/2}p_{j_n}(\lambda_{j_n}x).
\end{equation}
Applying Lemma~\ref{lem:moments} separately at orders $r$ and $s$ shows
that distinct tests have zero covariance under both measures. Indeed,
the tests are real, and for $k\ne l$ and $t\in\{r,s\}$,
\[
  \int\phi_k\phi_l\,d\mu^{*t}
  =\left(\int\phi_k\,d\mu^{*t}\right)
   \left(\int\phi_l\,d\mu^{*t}\right).
\]
For the diagonal terms, \eqref{eq:ab-ind} and
\eqref{eq:mean-r-ind} give
\begin{align}
  \int|\delta_j^{-1/2}p_j(\lambda_jx)|^2\,d\mu^{*r}(x)
  &\le\delta_j^{-1}\bigl((1-\delta_j)a_j+\delta_jb_j\bigr)\notag\\
  &\le\delta_j+1\le2.\label{eq:second-r-ind}
\end{align}
Similarly, by \eqref{eq:c-ind} and \eqref{eq:mean-s-ind},
\begin{align}
  \int|\delta_j^{-1/2}p_j(\lambda_jx)|^2\,d\mu^{*s}(x)
  &\le\delta_j^{-1}\bigl((1-\varepsilon_j)a_j+\varepsilon_jc_j\bigr)
       \notag\\
  &\le\delta_j+\frac{\varepsilon_j}{\delta_j}.
       \label{eq:second-s-ind}
\end{align}
The last bound is uniform on $J_r'$. Thus the covariance condition
\eqref{eq:Browncov} holds for the pair $(\mu^{*r},\mu^{*s})$.

Finally, \eqref{eq:mean-lower-ind} and \eqref{eq:phi-ind} imply
\[
  \left|\int\phi_n\,d\mu^{*r}-\int\phi_n\,d\mu^{*s}\right|
  \ge\frac12\sqrt{\delta_{j_n}}.
\]
The sum of squared differences diverges by \eqref{eq:deltadiv-ind}.
Brown's criterion therefore gives \mbox{$\mu^{*r}\perp\mu^{*s}$}.
\end{proof}

Independent powers alone give the remaining singularity assertion. The
following observation uses neither commutativity nor a Fourier decay
assumption.

\begin{lemma}\label{lem:independent-singular}
Let $G$ be a compact group with normalized Haar measure $m_G$, and let
$\nu$ be a probability measure on $G$. If
\[
  \nu^{*r}\perp\nu^{*s}\qquad(r,s\ge1,\ r\ne s),
\]
then $\nu^{*r}\perp m_G$ for every integer $r\ge1$.
\end{lemma}

\begin{proof}
Suppose that, for some $r\ge1$, the absolutely continuous part of
$\nu^{*r}$ is $f\,m_G$, where $f\ge0$ and $\int_G f\,dm_G>0$.
Set
\[
  \eta=\min\{f,1\}\,m_G,\qquad a=\eta(G)>0.
\]
Then $0\le\eta\le\nu^{*r}$ and $\eta\le m_G$. For $n\ge0$, put
$\eta_n=\eta*\nu^{*n}$, with $\nu^{*0}=\delta_e$, where $e$ is the
identity of $G$. Positivity of convolution and invariance of Haar measure
give
\[
  0\le\eta_n\le\nu^{*(r+n)},\qquad
  \eta_n\le m_G*\nu^{*n}=m_G,\qquad \eta_n(G)=a.
\]
The measures $\eta_n$ are pairwise mutually singular, since the
corresponding powers of $\nu$ are. Their densities with respect to $m_G$
are bounded by one and have pairwise disjoint sets of positivity, up to
null sets. Hence, for every $N\ge0$,
\[
  \sum_{n=0}^N\eta_n\le m_G,
  \qquad (N+1)a\le m_G(G)=1.
\]
This is impossible for sufficiently large $N$. Thus no positive
convolution power has a nonzero absolutely continuous part.
\end{proof}

Applying Lemma~\ref{lem:independent-singular} to the probability measure
$\mu$ from Proposition~\ref{prop:independent} gives
\begin{equation}\label{eq:allsingular}
  \mu^{*r}\perp m\qquad(r\ge1).
\end{equation}
Together with \eqref{eq:muMb}, Lemma~\ref{lem:Rajchman}, and
Proposition~\ref{prop:independent}, this completes the proof of
Theorem~\ref{thm:main}.

\section{Algebraic and spectral consequences}\label{sec:spectral}
We finish by describing the algebra generated by $\mu$, its spectrum in
$M(\T)$, and the associated convolution operator. The first consequence
is the standard isometric description of the algebra generated by an
independent-power measure; see \cite[p.~138]{GrahamMcGehee}.
We include the proof in the unital probability-measure setting.

\begin{corollary}[An isometric copy of $\ell^1(\mathbb N_0)$]
\label{cor:l1}
Set $\mu^{*0}=\delta_0$. The map
\[
  \Psi:\ell^1(\mathbb N_0)\longrightarrow M(\T),
  \qquad \Psi((a_n)_{n\ge0})=\sum_{n=0}^\infty a_n\mu^{*n},
\]
is an isometric unital Banach algebra homomorphism, where multiplication
on $\ell^1(\mathbb N_0)$ is convolution of sequences. Its image is the
closed unital subalgebra of $M(\T)$ generated by $\mu$.
\end{corollary}

\begin{proof}
Every positive power of $\mu$ belongs to $M_0(\T)$ and hence has no
atoms, by Wiener's theorem relating atoms to the mean squares of the
Fourier--Stieltjes coefficients; see \cite[Theorem~A.2.1, p.~415]{GrahamMcGehee}. In particular,
$\delta_0\perp\mu^{*n}$ for $n\ge1$. Proposition~\ref{prop:independent}
therefore shows that
\[
  \delta_0,\mu,\mu^{*2},\ldots
\]
are pairwise mutually singular positive probability measures. A finite
family of such measures admits pairwise disjoint Borel sets on which the
respective measures are concentrated. Consequently, for finite scalar
sequences,
\begin{equation}\label{eq:l1norm}
  \left\|\sum_{n=0}^N a_n\mu^{*n}\right\|_{M(\T)}
  =\sum_{n=0}^N|a_n|.
\end{equation}
The defining series for $\Psi$ converges absolutely in $M(\T)$, so
passage to the limit proves the isometry on all of
$\ell^1(\mathbb N_0)$. Multiplicativity follows from
$\mu^{*n}*\mu^{*m}=\mu^{*(n+m)}$ and absolute convergence. The image is
closed by isometry and is the closure of the polynomials in $\mu$,
since finitely supported sequences are dense in $\ell^1(\mathbb N_0)$.
\end{proof}

\begin{remark}[Two elementary $M_b$ consequences]\label{rem:Mb-algebra}
The space $M_b(\T)$ is a convolution ideal in $M(\T)$: for
$\nu\in M_b(\T)$ and $\eta\in M(\T)$,
\[
  S_N(\nu*\eta)=(S_N\nu)*\eta,
  \qquad \|\nu*\eta\|_b\le\|\nu\|_b\|\eta\|.
\]
Hence $\sup_{r\ge1}\|\mu^{*r}\|_b\le\|\mu\|_b$. Combined with
Corollary~\ref{cor:l1}, this shows that the closed nonunital algebra
spanned by $\{\mu^{*r}:r\ge1\}$ is a copy of $\ell^1(\mathbb N)$
contained in $M_b(\T)$. Indeed, $\|\nu\|\le\|\nu\|_b$ for
$\nu\in M_b(\T)$, by weak-$*$ convergence of $S_N\nu\,m$ to $\nu$.
More generally, if
$f(z)=\sum_{r\ge1}a_rz^r$ with $\sum_r|a_r|<\infty$, then
$f(\mu)=\sum_{r\ge1}a_r\mu^{*r}$ belongs to $M_b(\T)$. We record
these only as immediate by-products of the ideal property.
\end{remark}

\begin{remark}[The operator algebra]
The convolution representation $\nu\mapsto T_\nu$ is an isometric
unital homomorphism from $M(\T)$ into $B(L^1(\T))$. For completeness,
Young's inequality gives $\|T_\nu\|\le\|\nu\|$. Conversely, the Fej\'er
kernels satisfy $\|F_N\|_1=1$ and
$(\nu*F_N)m\to\nu$ weak-$*$, whence
\[
  \|\nu\|\le\liminf_N\|\nu*F_N\|_1\le\|T_\nu\|.
\]
Thus \eqref{eq:l1norm} also gives
\[
  \left\|\sum_{n=0}^N a_nT_\mu^n\right\|_{B(L^1(\T))}
  =\sum_{n=0}^N|a_n|.
\]
The closed unital operator algebra generated by $T_\mu$ is therefore
isometrically isomorphic to $\ell^1(\mathbb N_0)$.
\end{remark}

We next determine the spectrum in the full measure algebra. A classical
theorem states that a continuous probability measure on $\T$ with
independent powers has spectrum equal to the closed unit disc; see
\cite[Theorem~6.1.1(iv)]{GrahamMcGehee}. The circle-group case is
attributed there to Taylor; Bailey, Brown, and Moran \cite{BaileyBrownMoran}
developed the general spectral theory of independent-power measures.
Our measure is continuous because it belongs to $M_0(\T)$, and
Proposition~\ref{prop:independent} gives independent powers. Hence
\begin{equation}\label{eq:fulldisc}
  \sigma_{M(\T)}(\mu)=\{z\in\mathbb C:|z|\le1\}.
\end{equation}
Since $\widehat\mu(n)\to0$, the set
$\overline{\widehat\mu(\Z)}=\widehat\mu(\Z)\cup\{0\}$ is countable
and cannot equal this disc. This proves Corollary~\ref{cor:zafran}.
The full-disc assertion uses the cited theorem: the isometric
identification of the generated subalgebra with $\ell^1(\mathbb N_0)$
alone does not identify the spectrum in the ambient algebra $M(\T)$.

There is also a direct route to the assertion that the measure has non-natural spectrum,
which does not use the full-disc theorem. Zafran proved that if $G$ is a compact
abelian group, $0\ne\nu\in M_0(G)$, and every positive convolution power
of $\nu$ is singular with respect to Haar measure, then
\[
  \sigma_{M(G)}(\nu)\ne\widehat\nu(\widehat G)\cup\{0\};
\]
see \cite[Theorem~3.6]{Zafran}. Our measure satisfies these hypotheses
by \eqref{eq:allsingular}.

The obstruction can also be expressed directly in the quotient algebra.
Let $q:M(\T)\to M(\T)/L^1(\T)$ be the quotient map. Each $\mu^{*r}$
is a singular probability measure, so for every $f\in L^1(\T)$,
\[
  \|\mu^{*r}-f\,m\|=1+\|f\|_1.
\]
Taking the infimum over $f$ yields
\begin{equation}\label{eq:quotientnorm}
  \|q(\mu)^r\|
  =\dist_{M(\T)}(\mu^{*r},L^1(\T))=1\qquad(r\ge1).
\end{equation}
The spectral radius formula yields the remaining assertion in
\eqref{eq:introquotient}:
\begin{equation}\label{eq:quotientradius}
  r_{M(\T)/L^1(\T)}(q(\mu))=1.
\end{equation}

To relate this computation to Zafran's ideal, let $\Phi_{M(\T)}$ be the
maximal ideal space of $M(\T)$ and identify $n\in\Z$ with the character
$\nu\mapsto\widehat\nu(n)$. For
$\check\nu(\varphi)=\varphi(\nu)$, put
\[
  M_{00}(\T)=
  \{\nu\in M(\T):\check\nu(\varphi)=0
       \text{ for every }\varphi\in\Phi_{M(\T)}\setminus\Z\}.
\]
The characters outside $\Z$ are precisely those that annihilate
$L^1(\T)$. Hence $M_{00}(\T)$ is the inverse image under $q$ of the
Jacobson radical of $M(\T)/L^1(\T)$, or equivalently of its
quasinilpotent elements. In the notation fixed in the introduction,
\[
  M_{00}(\T)=\operatorname{Rad}(L^1(\T)).
\]
On the circle one has $M_{00}(\T)\subset M_0(\T)$, and Zafran's theorem
therefore identifies
\[
  \Ccal=M_{00}(\T)=\operatorname{Rad}(L^1(\T));
\]
see \cite{Zafran,HatoriSato}. Thus \eqref{eq:quotientradius} gives both
$\mu\notin\operatorname{Rad}(L^1(\T))$ and $\mu\notin\Ccal$.

Thus the construction gives, in particular,
$M_b(\T)\setminus\operatorname{Rad}(L^1(\T))\ne\varnothing$, answering
Andersson's stated nonemptiness question; the stronger identity
\eqref{eq:quotientnorm} also rules out finite convolution regularization.

The non-natural spectrum also has an operator-theoretic consequence.
Recall that a bounded operator $T$ on a complex Banach space $X$ is
\emph{decomposable} if, for every finite open cover
$\mathbb C=U_1\cup\cdots\cup U_k$, there are closed $T$-invariant
subspaces $X_1,\ldots,X_k$ such that
\[
  X=X_1+\cdots+X_k,
  \qquad \sigma(T|_{X_i})\subset U_i\quad(1\le i\le k).
\]
The sum here need not be direct. This is decomposability in local
spectral theory.

For a compact abelian group $G$ and a measure $\nu\in M_0(G)$,
Laursen and Neumann proved the equivalence
\[
  T_\nu\text{ is decomposable on }L^1(G)
  \quad\Longleftrightarrow\quad
  \nu\text{ has natural spectrum};
\]
see \cite[Theorem~4.5]{LaursenNeumann} and \cite{HatoriSato}. Applying
this result to our $\mu\in M_0(\T)$, whose spectrum is non-natural,
proves Corollary~\ref{cor:decomp}. The restriction to $M_0(G)$ is
important: natural spectrum does not characterize decomposability of
convolution operators throughout $M(G)$, as the examples in
\cite{HatoriSato} show.

Thus uniformly bounded Fourier partial sums are compatible not only
with persistent singularity and independent convolution powers, but
also with spectrum equal to the closed unit disc and failure of decomposability. The
uniform partial-sum bound imposes no convolution regularization of
the types considered above.

\section*{Use of artificial intelligence}
This work was supported by OpenAI through the ChatGPT for Academic
Researchers program. AI-assisted tools were used during the preparation
of this manuscript for mathematical verification and editorial support.
The author has independently verified the content and takes full
responsibility for its accuracy and conclusions.

\end{document}